\documentclass[12pt]{article}
\usepackage{amsmath,amsthm,amssymb,mathrsfs,color,bm}
\usepackage{mathrsfs}
\usepackage{iftex}
\ifpdf
  \usepackage{graphicx}% pdfLaTeX (arXiv)
\else
  \usepackage[dvipdfmx]{graphicx}% LaTeX + dvipdfmx
\fi
\usepackage{color}
\usepackage{here}
\usepackage{tikz}
\usetikzlibrary{patterns}
\usepackage{hyperref} % リンクを有効にするために必要
\IfFileExists{orcidlink.sty}{%
  \usepackage{orcidlink}%
}{%
  \definecolor{orcidlogocol}{RGB}{166,206,57}%
  \DeclareRobustCommand{\orcidlink}[1]{%
    \texorpdfstring{%
      \href{https://orcid.org/##1}{%
        \raisebox{-0.10ex}{%
          \textcolor{orcidlogocol}{%
            \textcircled{\raisebox{0.10ex}{\scriptsize\sffamily\bfseries iD}}%
          }%
        }%
      }%
    }{ORCID~##1}%
  }%
}
\newtheorem{tm}{Theorem}
[section]
\newtheorem{lm}[tm]{Lemma}

\newtheorem{exm}[tm]{Example}

\allowdisplaybreaks[4]
\makeatletter
\newcommand{\subscripts}[3]{%
  \@mathmeasure\z@\displaystyle{#2}%
  \global\setbox\@ne\vbox to\ht\z@{}\dp\@ne\dp\z@
  \setbox\tw@\box\@ne
  \@mathmeasure4\displaystyle{\copy\tw@_{#1}}%
  \@mathmeasure6\displaystyle{{#2}_{#3}}%
  \dimen@-\wd6 \advance\dimen@\wd4 \advance\dimen@\wd\z@
  \hbox to\dimen@{}\mathop{\kern-\dimen@\box4\box6}%
}
\makeatother
\newcommand{\III}{{\vert \kern-.10em \vert \kern-.10em \vert}}

\makeatletter
 
 \@addtoreset{equation}{section}
\makeatother
\newcommand{\tp}[1]{{}^{t}\!#1}
\begin{document}
%%%%%%%%%%%%%%%%%%%%%%%%%%%%%%%%%%%%%%%%%%%
\setlength{\baselineskip}
{15.5pt}
%{15.5pt}
%%%%%%%%%%%%%%%%%%%%%%%%%%%%%%%%%%%%%%%%%%%
\title{A generalization of L\'{e}vy's stochastic area formula \\
for Gaussian rough paths}
\author{
\large{
Hiroshi Kawabi\,\orcidlink{0009-0006-7128-206X}
\hspace{-2mm}
\thanks{Department of Mathematics, Hiyoshi Campus, Keio University, 4-1-1 Hiyoshi, Kohoku-ku, 
Yokohama, Kanagawa 223-8521,
Japan (e-mail address: {\tt kawabi@keio.jp})}~
}}
%%%%%%%%%%%%%%%%%%%%%%%%
\date{\today}
\maketitle 
%\vspace{-5mm}
\begin{abstract}
In this paper, we establish a generalization of L\'{e}vy's stochastic
area formula for the canonical rough path lift of a \(2n\)-dimensional
centered Gaussian process with independent, identically distributed
components. The formula extends the identity obtained by Aihara, Akahori,
Fujii and Nitta \cite{AAFN} for Brownian motion. We compute the expectation
of the exponential of a functional involving stochastic areas and quadratic
forms in the terminal values. The result is expressed in terms of spectral
data of a Hilbert--Schmidt operator associated with the
stochastic area and an \(n\times n\) determinant. We also give a necessary
and sufficient condition for absolute integrability in terms of the
covariance at the terminal time and the symmetric part of the quadratic form
in the terminal values.  As examples, we discuss 
the cases of fractional Brownian motion with Hurst parameter \(H\in(1/4,1)\), Brownian
motion and the Brownian bridge. In particular, we recover the
formula in \cite{AAFN} under the exact condition for absolute integrability.
\vspace{2mm} 
\vspace{1mm} 
\end{abstract}
%%%%%%%%%%%%%%%%%%%%%%%%%%%%%%%%%%%%%%%%%%%%
\section{Introduction}
\label{sec:Introduction}
%%%%%%%%%%%%%%%%%%%%%%%%%%%%%%%%%%%%%%%%%%%%
Let \(B=(B^{1}_{t},B^{2}_{t})_{0\leq t\leq 1}\) be a two-dimensional Brownian motion starting at
the origin. The stochastic area
\[
S:=\frac{1}{2}\int_{0}^{1}(B^{1}_{t}\,dB^{2}_{t}-B^{2}_{t}\,dB^{1}_{t})
\]
enclosed by the path of \(B\) and its chord was introduced by Paul L\'evy
\cite{Le40}. Later, in \cite{Le51}, he obtained the celebrated formula
\[
\mathbb{E}[e^{\sqrt{-1}\xi S}]=\Bigl(\cosh\frac{\xi}{2}\Bigr)^{-1},
\qquad \xi\in\mathbb{R},
\]
and also its conditional version given the terminal value \(B_1\).
Since then, {\emph{L\'evy's stochastic area formula}} has been studied and
extended by many authors. For example, Yor \cite{Yo} gave another proof,
Gaveau \cite{Ga} related it to the heat kernel on the Heisenberg group,
Helmes and Schwane \cite{HS} and Geng and Qian \cite{GQ} studied the higher-dimensional case, and
Ikeda, Kusuoka and Manabe \cite{IKM} extended it to Gaussian processes.

On the other hand, {\emph{quadratic Wiener functionals}} also appear in soliton
theory. It is known that reflectionless potentials and soliton solutions
of the KdV equation are represented by quadratic Wiener functionals. See
Ikeda and Taniguchi \cite{IT} and the references therein. In this
direction, Aihara, Akahori, Fujii and Nitta \cite{AAFN} generalized
L\'evy's stochastic area formula in the following way. For \(n\) independent
two-dimensional Brownian motions, they considered a quadratic Wiener functional
consisting of stochastic areas and quadratic forms in the terminal
values coupled through a real \(n\times n\) matrix, and proved that the
expectation of the exponential of this functional is equal to the
reciprocal of an explicit \(n\times n\) determinant. 
In this paper, we call this identity the {\emph{AAFN formula}} (see 
Example~\ref{exam:bm} below for the precise form).
They also showed that this formula gives representations of the
tau functions of KP solitons.

The purpose of this paper is to study this kind of formula for {\emph{Gaussian
rough paths}}. For Gaussian processes other than
Brownian motion, stochastic areas can be defined in the framework of
rough path theory. Coutin and Qian \cite{CQ} constructed the rough path
lift of fractional Brownian motion with Hurst parameter \(H>1/4\). 
In \cite{FV}, Friz and Victoir  extended this construction to Gaussian processes
whose covariance has finite {\emph{2D $\rho$-variation}}
with \(1\leq \rho<2\) (see also \cite{FV10, FR}). Ferreiro-Castilla and Utzet \cite{FCU} represented the
stochastic area of Gaussian processes as a double Wiener--It\^o integral and
studied its characteristic function. 

In this paper, we consider a \(2n\)-dimensional Gaussian process with
independent, identically distributed components, and take the stochastic
areas from its canonical rough path lift. The main point of our result
(Theorem~\ref{thm:main})
is to incorporate quadratic terminal terms coupled through an arbitrary
real matrix into the stochastic area formula for
Gaussian processes.
The resulting expression \eqref{eq:main} separates the
contribution of the stochastic areas from that of the terminal values.
The former is described by the spectral data of a Hilbert--Schmidt
operator, while the latter is expressed through an \(n\times n\) determinant
involving a scalar resolvent quantity associated with the terminal value.
The Wiener-chaos representation and the characteristic function of the
uncoupled stochastic area are covered by \cite{FCU}. 
We also identify the exact condition for absolute integrability.

In our proof of Theorem~\ref{thm:main}, we reduce the problem to a finite-dimensional complex
Gaussian integral. 
In evaluating the resulting determinant, {\emph{Sylvester's determinant
identity}} plays a crucial role.
In Section~\ref{sec:examples}, we discuss fractional
Brownian motion with \(H\in(1/4,1)\), Brownian motion and the Brownian
bridge as examples.
In particular, we recover the AAFN formula under the exact condition 
for absolute integrability, which is weaker than 
the smallness assumptions in \cite{AAFN}.

%%%%%%%%%%%%%%%%%%%%%%%%%%%%%%%%%
\section{Framework and the results}
\label{sec:framework}
%%%%%%%%%%%%%%%%%%%%%%%%%%%%%%%%%
%
For a matrix \(A\), we denote by
\(\tp{A}\) and \(A^{*}\) the transpose and the conjugate transpose of
\(A\), respectively. We also denote by \(A^{*}\) the adjoint of an
operator \(A\) on a Hilbert space.
We denote by \(I_m\) the \(m\times m\) identity matrix.
Vectors are regarded as column vectors. We use
\((x,y)_{\mathbb R^m}=\sum_{i=1}^m x^iy^i\) and
\((u,v)_{\mathbb C^m}=\sum_{i=1}^m u^i\overline{v^i}\)
for $x=\tp \,{(x^{1},\ldots, x^{m})}, y=\tp \,{(y^{1},\ldots, y^{m})}\in {\mathbb R}^{m}$
and 
$u=\tp \,{(u^{1},\ldots, u^{m})}, v=\tp \,{(v^{1},\ldots, v^{m})}\in {\mathbb C}^{m}$.
For a Hermitian matrix \(A\), we write \(A>0\) if \(A\) is positive definite.

We fix \(n\in\mathbb N\) and 
let \(w=(w_t)_{0\leq t\leq1}\) with
\(
w_t=\tp \,{(\beta^{1}_t,\ldots,\beta^{n}_t,
\widetilde{\beta}^{1}_t,\ldots,\widetilde{\beta}^{n}_t)}\), \(0\leq t \leq 1\),
be a \(2n\)-dimensional centered
continuous Gaussian process starting at zero with i.i.d.\ components,
defined on a complete probability space
$(\Omega, {\cal F}, \mathbb P)$. 
We write \(\beta_t= \tp \,(\beta^{1}_t,\ldots,\beta^{n}_t)\) and 
\(\widetilde{\beta}_t=\tp \,(\widetilde{\beta}^{1}_t,\ldots,\widetilde{\beta}^{n}_t)\).
We assume that 
the covariance function 
$$R(s,t):={\mathbb E}[\beta^{1}_{s}\beta^{1}_{t}], \quad (s,t)\in [0,1]^{2}$$
has finite {\emph{2D \(\rho\)-variation}} with some $1\leq \rho <2$, namely, 
\[
V_\rho(R;[0,1]^2)
:=\sup_{\mathcal D,\mathcal D'}
\Big(\sum_{i,j}|(\delta R)_{ij}|^\rho \Big)^{1/\rho}<\infty
\]
for some $1\leq \rho <2$, 
where \(\mathcal D=(s_i)\) and \(\mathcal D'=(t_j)\) are finite
partitions of \([0,1]\) and
\[
(\delta R)_{ij}:=R(s_i,t_j)-R(s_{i-1},t_j)-R(s_i,t_{j-1})+R(s_{i-1},t_{j-1}).
\]
We put
\[
v:=R(1,1)\geq0.
\]

We now take \(p\in(2\rho,4)\) and denote by
\(\lfloor p\rfloor\) the greatest integer not exceeding
\(p\). 
Let \(\mathbf{W}=(1,{W}^{1},\ldots,{W}^{\lfloor p\rfloor})\) 
be the Gaussian rough path, that is, the 
canonical geometric \(p\)-rough path lift of the Gaussian process \(w\), which is constructed as
the limit of the dyadic piecewise linear approximations. It is known
that these approximations converge in \(L^{q}(\mathbb P)\) with respect to the
\(p\)-variation rough path topology for every \(q\geq 1\)
(see, e.g., \cite[Theorem 15.33]{FV10} and \cite[Theorem 1]{FR}). 
Note that this rough path lift has a third level path when $3/2 \leq \rho<2$, and that
\({W}^{1}_{0,t}=w_{t}=\tp \,(\tp \, \beta_t, \tp \, \widetilde{\beta}_t), \, 0\leq t\leq 1\).
%%%%%%%%%%%%%%%%%
We define 
\begin{equation}\label{eq:area}
\begin{aligned}
\mathcal{A}^{l}(t)
:=
{W}^{2;n+l,l}_{0,t}-{W}^{2;l,n+l}_{0,t}, \quad
S^{l}:=-\tfrac{1}{2}\mathcal{A}^{l}(1), \qquad 1\leq l \leq n,
\end{aligned}
\end{equation}
%%%%%%%%%%%%%%%%%%%%%%%%%%%%%
where 
\({W}^{2;i,j}\) stands for the \((i,j)\)-component
of the second level path \({W}^{2}\).

Let \(\mathcal{H}=\mathcal{H}^R\) be the Hilbert space obtained as the completion of
the space of step functions on \([0,1]\) modulo null vectors with
respect to the inner product determined by
\[
(\mathbf{1}_{[0,s]},\mathbf{1}_{[0,t]})_{\mathcal{H}}=R(s,t).
\]
Note that \(\mathcal H\) is separable and
\({\rm{dim}}\, \mathcal{H}\) may be finite.
We write \(d:=\dim\mathcal H\in {\mathbb N} \cup\{\infty\}\), and
set \(e:=\mathbf{1}_{[0,1]}\in\mathcal{H}\).
Then we have
\(\|e\|^{2}_{\mathcal{H}}=v\). Let \(X_l\) and \(Y_l\),
\(1\leq l\leq n\), be independent isonormal Gaussian processes over
\(\mathcal{H}\) such that
\[
\beta^{l}_{t}=X_l(\mathbf{1}_{[0,t]}),
\qquad
\widetilde{\beta}^{l}_{t}=Y_l(\mathbf{1}_{[0,t]}).
\]
Then we have the following lemma.
%%%%%%%%%%%%%%%%%%%%%%%%%%%%%%%%%%%%%%%
\begin{lm}
\label{lem:area-operator}
There exists a unique real skew-symmetric Hilbert--Schmidt operator
\(K\) on \(\mathcal{H}\), that is, \(K^{*}=-K\), such that, for \(1\leq l\leq n\),
\begin{equation}\label{eq:area-series}
S^{l}=\sum_{j,k=1}^{d} \big(K\varphi_k,\varphi_j\big)_{\mathcal{H}}
X_l(\varphi_j)Y_l(\varphi_k)
\end{equation}
holds in \(L^{2}(\mathbb P)\) for every real complete orthonormal basis \((\varphi_j)_{j=1}^{d}\) of
\(\mathcal{H}\), where
the double sum on the right-hand side
denotes the \(L^2(\mathbb P)\)-limit of the partial sums over \(1\leq j,k\leq N\)
as \(N\to\infty\) in the case \(d=\infty\).
Moreover, we have
\begin{equation}
\label{eq:HS}
\|K\|^{2}_{\rm{HS}(\cal H)}=\mathbb{E}[(S^{l})^{2}].
\end{equation}
\end{lm}
%%%%%%%%%%%%%%%%%%%%%%%%%%%%%%%%%%%%%%
Since \(K\) is real and skew-symmetric, the nonzero eigenvalues of
its complexification
appear in pairs \({\sqrt{-1}}\mu_k,-{\sqrt{-1}}\mu_k\) with
\(\mu_{k}>0\).
We denote by \(r_K\) the number of these
pairs, counted with multiplicity.
Then we see that
\(r_K=\tfrac12\operatorname{rank}K\) if \(K\) has finite rank, and
\(2 \sum_{k=1}^{r_K}\mu_k^{2}=\|K\|^{2}_{\mathrm{HS}(\cal H)}<\infty\). 

For
\(a\in\mathbb{R}\), we define
\begin{equation}\label{eq:spectral-data}
\begin{aligned}
\Delta_R(a)&:=\prod_{k=1}^{r_K}(1+a^{2}\mu_k^{2})
=\det(I_{\mathcal H}+a^{2}K^{*}K)^{1/2},
\\
m_R(a)&:=\big((I_{\mathcal H}+a^{2}K^{*}K)^{-1}e,e
\big)_{\mathcal{H}}.
\end{aligned}
\end{equation}
Here the Fredholm determinant is well-defined because \(K^{*}K\) is
of trace class, and we take its positive square root. The product in
\eqref{eq:spectral-data} may be a finite product or the empty product.
Both functions are even, and we have
\[
\Delta_R(0)=1,\qquad m_R(0)=v,\qquad 0\leq m_R(a)\leq v.
\]

Let \(\Lambda=\operatorname{diag}(\lambda_1,\ldots,\lambda_n)\) with
\(\lambda_l>0\), and let \(C\) be a real \(n\times n\) matrix. We define
\begin{align*}
M&=M_R(\Lambda):=\operatorname{diag}(m_R(2\lambda_1),\ldots,m_R(2\lambda_n)),
\nonumber \\
D_C&:=\Lambda^{1/2}\frac{C-\tp \,{C}}{2}\Lambda^{1/2}, \quad
E_C:=\Lambda^{1/2}\frac{C+\tp \,{C}}{2}\Lambda^{1/2}.
\end{align*}
Although these matrices depend on \(\Lambda\), we do not indicate this dependence
for notational simplicity. We consider the complex-valued functional
\begin{equation}\label{eq:functional}
\begin{aligned}
\mathcal{S}^{R}_{C,\Lambda}
:={}&\sqrt{-1}\Big\{\sum_{l=1}^{n}\lambda_l\mathcal{A}^{l}(1)
+\big(D_C\beta_1,\widetilde{\beta}_1\big)_{\mathbb R^n}\Big\}\\
&+\frac{\, 1\, }{2}\Big \{
\big(E_C\beta_1,{\beta}_1\big)_{\mathbb R^n}
+
\big(E_C{\widetilde{\beta}}_1,\widetilde{\beta}_1\big)_{\mathbb R^n}
\Big\}.
\end{aligned}
\end{equation}

We are now in a position to state the main result of this paper.
%%%%%%%%%%%%%%%%%%%%%%%%%%%%%%%%%%%%%
\begin{tm}
[A generalized AAFN formula for Gaussian rough paths]
\label{thm:main}
The random variable \(e^{\mathcal{S}^{R}_{C,\Lambda}}\) is absolutely
integrable if and only if \(I_n-vE_C>0\). Under this condition, we have
\begin{equation}\label{eq:main}
\mathbb{E}\big[ e^{\mathcal{S}^{R}_{C,\Lambda}} \big]
=
\frac{
\prod_{l=1}^{n}\Delta_R(2\lambda_l)^{-1}
}
{\det(I_n-C\Lambda M)}.
\end{equation}
Here the denominator is strictly positive. Moreover, we have
\[
\mathbb{E}\big[ \vert e^{\mathcal{S}^{R}_{C,\Lambda}}\vert \big]=\det(I_n-vE_C)^{-1}.
\]
In particular, if \(v=0\), then \(e\) is the zero element of
\(\mathcal H\) and \(m_R=0\). Thus \(M=0\), and
\eqref{eq:main} holds for every real \(n\times n\) matrix \(C\),
with the denominator equal to \(1\).
\end{tm}

Before closing this section, let us consider the case \(C=0\) and \(n=1\). In this case, the condition
\(I_{n}-vE_C=1>0\) is always satisfied. Since
\(\mathcal{A}^{1}(1)=-2S^{1}\) and \(\Delta_R\) is even,
\eqref{eq:main} reads
\[
\mathbb{E}[e^{\sqrt{-1}aS^{1}}]=\Delta_R(a)^{-1},\qquad a\in\mathbb{R}.
\]
This formula corresponds to L\'evy's stochastic area formula 
in \cite[Proposition 4.1]{FCU}
for the two-dimensional
Gaussian process
\(w=(\beta^{1}_{t},\widetilde{\beta}^{1}_{t})_{0\leq t \leq 1}\).
%
%%%%%%%%%%%%%%%%%%%%%%%%%%%%%%%%%%%%%%%%%%%%%%%
\section{Proof of the results}\label{sec:proofs}
%%%%%%%%%%%%%%%%%%%%%%%%%%%%%%%%%%%%%%%%%%%%%%%
In this section, we prove the results stated in
Section~\ref{sec:framework}. We begin with Lemma~\ref{lem:area-operator}.
\begin{proof}[\bfseries{Proof of Lemma~\ref{lem:area-operator}}]
For \(r\in\mathbb{N}\), we consider the dyadic partition
\(\{t_i:=i2^{-r}; 0\leq i\leq 2^{r} \}\) of \([0,1]\), and we put
\(h_i=\mathbf{1}_{[0,t_{i-1}]}\) and \(d_i=\mathbf{1}_{(t_{i-1},t_i]}\)
for \(1\leq i \leq2^r\).
Let \(S^{l,(r)}\) be the stochastic area of the \(r\)-th dyadic
piecewise linear approximation of the two-dimensional Gaussian process
\((\beta^{l},\widetilde{\beta}^{l})\).
With the normalization in \eqref{eq:area}, we have
\[
S^{l,(r)}
=\frac{1}{2}\sum_{i=1}^{2^r}\bigl\{X_l(h_i)Y_l(d_i)
-X_l(d_i)Y_l(h_i)\bigr\}.
\]
In view of this expression, we define a skew-symmetric operator \(L_r\)
of finite rank by
\[
L_r:=\frac{1}{2}\sum_{i=1}^{2^r}(h_i\otimes d_i-d_i\otimes h_i),
\]
where \((f\otimes g)u:=(u,g)_{\mathcal{H}}f\) for $f,g,u\in {\cal H}$.
Let \((\varphi_j)_{j=1}^{d}\) be a real complete orthonormal basis of
\(\mathcal{H}\). By expanding \(h_i\) and \(d_i\) with respect to this
basis, we have
\begin{equation}\label{eq:area-series-r}
\begin{aligned}
S^{l,(r)}
&=\frac{1}{2}\sum_{i=1}^{2^r}\sum_{j,k=1}^{d}
\bigl\{(\varphi_j,h_i)_{\mathcal{H}}(\varphi_k,d_i)_{\mathcal{H}}
-(\varphi_j,d_i)_{\mathcal{H}}(\varphi_k,h_i)_{\mathcal{H}}\bigr\}
X_l(\varphi_j)Y_l(\varphi_k)\\
&=\sum_{j,k=1}^{d}\big(\varphi_j,L_r \varphi_k\big)_{\mathcal{H}}
X_l(\varphi_j)Y_l(\varphi_k).
\end{aligned}
\end{equation}

On the other hand, since \(X_l\) and \(Y_l\) are independent, the family
\(\big(X_l(\varphi_j)Y_l(\varphi_k)\big)_{1\leq j,k\leq d}\) is an orthonormal
system in \(L^{2}(\mathbb{P})\). Hence, for every real array
\((c_{jk})_{1\leq j,k\leq d}\) with \(\sum_{j,k=1}^{d}|c_{jk}|^{2}<\infty\),
the series \(\sum_{j,k=1}^{d}c_{jk}X_l(\varphi_j)Y_l(\varphi_k)\) converges
in \(L^{2}(\mathbb{P})\), and we have
%%%%%
\begin{equation}\label{eq:isometry}
\mathbb{E}\bigg[ \, \Big|\sum_{j,k=1}^{d}c_{jk}X_l(\varphi_j)Y_l(\varphi_k)\Big|^{2} \, \bigg]
=\sum_{j,k=1}^{d}|c_{jk}|^{2}.
\end{equation}
%%%%%
In particular, the series in \eqref{eq:area-series-r} converges in
\(L^{2}(\mathbb{P})\). Indeed, since $L_r$ is of finite rank, it is a Hilbert--Schmidt operator, 
and hence $\sum_{j,k=1}^{d}(\varphi_j,L_r\varphi_k)_{\mathcal{H}}^{2}=\|L_r\|^{2}_{\mathrm{HS}(\cal H)}<\infty$.
Applying \eqref{eq:isometry} to
\(c_{jk}=(\varphi_j,(L_r-L_{r'})\varphi_k)_{\mathcal{H}}\) and using
\eqref{eq:area-series-r}, we obtain
\[
\|L_r-L_{r'}\|^{2}_{\mathrm{HS}(\cal H)}
=\mathbb{E}[\, |S^{l,(r)}-S^{l,(r')}|^{2}\, ],\quad r,r'\in\mathbb{N}.
\]
Since \(S^{l,(r)}\) converges to \(S^{l}\) in \(L^{2}(\mathbb P)\) as
\(r\to\infty\), \((L_r)_{r\in \mathbb N}\) is a Cauchy sequence with respect to the
Hilbert--Schmidt norm. We denote its limit by \(K\), which is also real and
skew-symmetric. Applying \eqref{eq:isometry} to
\(c_{jk}=(\varphi_j,(L_r-K)\varphi_k)_{\mathcal{H}}\), we see that the
right-hand side of \eqref{eq:area-series-r} converges to that of
\eqref{eq:area-series} in \(L^{2}(\mathbb P)\) as \(r\to\infty\). Since
the left-hand side of \eqref{eq:area-series-r} converges to \(S^{l}\),
we obtain \eqref{eq:area-series}. Finally, \eqref{eq:HS} and the
uniqueness of \(K\) follow from \eqref{eq:isometry}.
\end{proof}
%%%%%%%%%%%%%%%%%%%%%%%%%%%%%%%%%%%%%%%%%%%
Next, we prepare a lemma which will be used in the proof of Theorem~\ref{thm:main}.
%%%%%%%%%%%%%%%%%%%%%%%%%%%%%%%%%%%%%%%%%%%
\begin{lm}
\label{lem:complex-gauss}
Let \(q\in\mathbb N\) and let \(Z\) be a random vector in \(\mathbb C^q\)
whose components are independent standard complex normal random
variables. Namely, each component is of the form
\((X+\sqrt{-1}Y)/\sqrt{2}\) with independent \(N(0,1)\) random variables
\(X\) and \(Y\). Let \(L\) be a \(q\times q\) complex matrix, and write
\[
L=A+\sqrt{-1}B,\quad
A:=\frac{L+L^{*}}{2},\quad B:=\frac{L-L^{*}}{2\sqrt{-1}},
\]
where \(A\) and \(B\) are Hermitian matrices. If \(I_q-A>0\), 
then \(e^{(LZ,Z)_{\mathbb C^q}}\) is absolutely
integrable and
\begin{equation}\label{eq:complex-gauss}
\mathbb{E}[e^{(LZ,Z)_{\mathbb C^q}}]=\det(I_q-L)^{-1}.
\end{equation}
\end{lm}
%%%%%%%%%%%%%%%%%%%%%%%%%%%%%%%%%%%%%%%%%
\begin{proof}
We identify \(\mathbb C^{q}\) with \(\mathbb R^{2q}\) by writing
\(z=x+\sqrt{-1}\,y\) with \(x,y\in\mathbb R^{q}\). Then the law of
\(Z\) has the density \(\pi^{-q}e^{-|z|^{2}}\) with respect to the
Lebesgue measure \(dx\,dy\) on \(\mathbb R^{2q}\). We also recall that, for
\(c\in\mathbb C\) with \(\operatorname{Re}c>0\),
\begin{equation}\label{eq:gauss-1dim}
\frac{1}{\pi}\int_{\mathbb R^{2}}e^{-c(s^{2}+t^{2})}\,dsdt
=\frac{1}{c}.
\end{equation}
Since \(A\) and \(B\) are Hermitian, both \((AZ,Z)_{\mathbb C^q}\) and
\((BZ,Z)_{\mathbb C^q}\) are real-valued. Hence we have
\[
(LZ,Z)_{\mathbb C^q}
=(AZ,Z)_{\mathbb C^q}+\sqrt{-1}\,(BZ,Z)_{\mathbb C^q}, \quad
|e^{(LZ,Z)_{\mathbb C^q}}|
=e^{(AZ,Z)_{\mathbb C^q}}.\]

We put \(G:=I_q-A>0\). Then we have
\(|z|^{2}-(Az,z)_{\mathbb C^q}=(Gz,z)_{\mathbb C^q}\) for
\(z\in\mathbb C^{q}\). We consider the change of variables
\(z=G^{-1/2}u\), where \(u=\sigma+\sqrt{-1}\,\tau\) with
\(\sigma,\tau\in\mathbb R^{q}\). This is a linear transformation on
\(\mathbb R^{2q}\) whose Jacobian is
\(|\det G^{-1/2}|^{2}=(\det G)^{-1}\), and we have
\((Gz,z)_{\mathbb C^q}=|u|^{2}\). Hence we obtain
\begin{align*}
\mathbb{E}[\,|e^{(LZ,Z)_{\mathbb C^q}}|\,]
&=\frac{1}{\pi^{q}}\int_{\mathbb R^{2q}}
e^{-(Gz,z)_{\mathbb C^q}}\,dx dy\\
&=\frac{1}{\det G}\cdot\frac{1}{\pi^{q}}
\int_{\mathbb R^{2q}}e^{-|u|^{2}}\,d\sigma d\tau
=\frac{1}{\det G}<\infty.
\end{align*}
This shows the absolute integrability.

Next, we prove \eqref{eq:complex-gauss}. We put
\(\widetilde{B}:=G^{-1/2}BG^{-1/2}\). By the same change of variables,
we have
\begin{align*}
\mathbb{E}[e^{(LZ,Z)_{\mathbb C^q}}]
&=\frac{1}{\pi^{q}}\int_{\mathbb R^{2q}}
e^{-(Gz,z)_{\mathbb C^q}+\sqrt{-1}\,(Bz,z)_{\mathbb C^q}}\,dxdy\\
&=\frac{1}{\det G}\cdot\frac{1}{\pi^{q}}\int_{\mathbb R^{2q}}
e^{-|u|^{2}+\sqrt{-1}\,(\widetilde{B}u,u)_{\mathbb C^q}}\,d\sigma d\tau.
\end{align*}
Since \(\widetilde{B}\) is Hermitian, there exists a unitary matrix \(U\)
such that
\(U^{*}\widetilde{B}U=\operatorname{diag}(b_1,\ldots,b_q)\) with
\(b_1,\ldots,b_q\in\mathbb R\). We consider the change of variables
\(u=U\zeta\), where \(\zeta=s+\sqrt{-1}\,t\) with
\(s=\tp \, {(s_1,\ldots,s_q)},t=\tp \,{(t_1,\ldots,t_q)}\in\mathbb R^{q}\).
This is an orthogonal transformation on \(\mathbb R^{2q}\), and we have
\[
|u|^{2}=\sum_{j=1}^{q}(s_j^{2}+t_j^{2}),\qquad
(\widetilde{B}u,u)_{\mathbb C^q}
=(U^{*}\widetilde{B}U\zeta,\zeta)_{\mathbb C^q}
=\sum_{j=1}^{q}b_j(s_j^{2}+t_j^{2}).
\]
Hence, by Fubini's theorem and \eqref{eq:gauss-1dim} with
\(c=1-\sqrt{-1}\,b_j\), we obtain
\begin{align*}
\mathbb{E}[e^{(LZ,Z)_{\mathbb C^q}}]
&=\frac{1}{\det G}\prod_{j=1}^{q}\frac{1}{\pi}\int_{\mathbb R^{2}}
e^{-(1-\sqrt{-1}\,b_j)(s_j^{2}+t_j^{2})}\,ds_j dt_j
=\frac{1}{\det G}\prod_{j=1}^{q}\frac{1}{1-\sqrt{-1}\,b_j}.
\end{align*}
On the other hand, we have
\begin{align*}
\det(I_q-L)
&=\det(G-\sqrt{-1}\,B)\\
&=\det\bigl\{G^{1/2}(I_q-\sqrt{-1}\,\widetilde{B})G^{1/2}\bigr\}
=\det G\prod_{j=1}^{q}(1-\sqrt{-1}\,b_j).
\end{align*}
By combining these identities, we obtain \eqref{eq:complex-gauss}.
\end{proof}
%%%%%%%%%%%%%%%%%%%%%%%%%%%%%%%%%%%
We now prove the main result.

\begin{proof}[{\bfseries Proof of Theorem~\ref{thm:main}}]
We first prove the assertions concerning absolute integrability.
By \eqref{eq:functional}, we have
\[
|e^{\mathcal{S}^{R}_{C,\Lambda}}|
=\exp\Big[\frac12\bigl\{(E_C\beta_1,\beta_1)_{\mathbb R^n}
+(E_C\widetilde\beta_1,\widetilde\beta_1)_{\mathbb R^n}\bigr\}\Big].
\]
Let \(\kappa_1,\ldots,\kappa_n\) be the eigenvalues of the real
symmetric matrix \(E_C\), and take an orthogonal matrix \(O\) such that
\(\tp \,{O}E_CO=\operatorname{diag}(\kappa_1,\ldots,\kappa_n)\).
Since \(\beta_1\) and \(\widetilde{\beta}_1\) are independent Gaussian
random vectors with law \(N(0,vI_n)\), so are \(\tp \,{O}\beta_1\) and
\(\tp \,{O}\widetilde{\beta}_1\). Hence we have
\[
\mathbb{E}[\,|e^{\mathcal{S}^{R}_{C,\Lambda}}|\,]
=\prod_{i=1}^{n}\Bigl(\mathbb{E}\bigl[e^{v\kappa_i g^{2}/2}\bigr]\Bigr)^{2},
\]
where \(g\) is an \(N(0,1)\) random variable. Recalling
\[
\mathbb{E}\bigl[e^{\theta g^{2}/2}\bigr]
=\begin{cases}
(1-\theta)^{-1/2}, & \theta<1,\\
\infty, & \theta\geq 1,
\end{cases}
\]
we see that \(e^{\mathcal{S}^{R}_{C,\Lambda}}\) is absolutely integrable
if and only if \(v\kappa_i<1\) for all \(1\leq i\leq n\), that is,
\(I_n-vE_C>0\). Under this condition, we obtain
\[
\mathbb{E}[\,|e^{\mathcal{S}^{R}_{C,\Lambda}}|\,]
=\prod_{i=1}^{n}(1-v\kappa_i)^{-1}
=\det(I_n-vE_C)^{-1}.
\]

Next, we prove \eqref{eq:main}. If \(d<\infty\), the finite-dimensional
approximation below is not necessary, because the following argument
with \(N=d\) gives \eqref{eq:main} directly. 
Hence we only consider the case \(d=\infty\)
in what follows. 

We take a real complete orthonormal basis
\((\varphi_j)_{j=1}^{\infty}\) of \(\mathcal{H}\)
such that \(e=\sqrt v\,\varphi_1\). When \(v=0\), any real complete
orthonormal basis satisfies this condition.
For \(N\in\mathbb N\), we put
\(\mathcal H_N:=\operatorname{Span}\{\varphi_1,\ldots,\varphi_N\}\),
and let \(\pi_N\) be the orthogonal projection onto \(\mathcal H_N\).
Then we have \(\pi_N e=e\). We regard
\(K_N:=\pi_NK\pi_N\) as a real skew-symmetric \(N\times N\) matrix
on \(\mathcal H_N\) using this basis. When taking operator
norms or limits on \(\mathcal H\), we extend this matrix by zero on
\(\mathcal H_N^\perp\). 
We set
\[
S^{l}_{N}:=\sum_{j,k=1}^{N}\big(K\varphi_k,\varphi_j\big)_{\mathcal{H}}
X_l(\varphi_j)Y_l(\varphi_k), \quad 1\leq l\leq n,
\]
and define
\begin{align*}
\mathcal{S}_N
:={}&\sqrt{-1}\Big\{-2\sum_{l=1}^{n}\lambda_l S^{l}_{N}
+\big(D_C\beta_1,\widetilde{\beta}_1\big)_{\mathbb R^n}\Big\}
+\frac{\, 1\, }{2}\Big \{
\big(E_C\beta_1,{\beta}_1\big)_{\mathbb R^n}
+
\big(E_C{\widetilde{\beta}}_1,\widetilde{\beta}_1\big)_{\mathbb R^n}
\Big\}.
\end{align*}
Namely, \(\mathcal{S}_N\) is obtained by replacing
\(\mathcal{A}^{l}(1)=-2S^{l}\) in \eqref{eq:functional} by
\(-2S^{l}_{N}\). We also set
\begin{align*}
Z_{lj}&:=\frac{X_l(\varphi_j)+\sqrt{-1}Y_l(\varphi_j)}{\sqrt{2}}, \quad
Z_{l}:=\tp \,{(Z_{l1},\ldots,Z_{lN})}\in\mathbb C^N, \quad
1\leq l \leq n,~1\leq j \leq N, \\
Z&:=\tp \,{(\tp \,{Z_1},\ldots,\tp \,{Z_n})}\in\mathbb C^{nN}, \quad 
P:=\sqrt{v}\,I_n\otimes\tp \,{\mathbf{e}_1}, \quad
B_C:=E_C-D_C=\Lambda^{1/2}\,(\tp \,{C})\Lambda^{1/2},
\end{align*}
where \(\mathbf{e}_1=\tp \,{(1,0,\ldots,0)}\in\mathbb{R}^{N}\) and
the tensor symbol \(\otimes\) denotes the Kronecker product of matrices.
In particular, \(P\in\mathbb R^{n\times nN}\).
Then we have
\[PZ=\frac{1}{\sqrt 2}(\beta_1+\sqrt{-1}\widetilde{\beta}_1).\] 
We also write
\(X_l\) and \(Y_l\) for the column vectors
\(\tp \,{(X_l(\varphi_j))_{j=1}^N}\) and
\(\tp \,{(Y_l(\varphi_j))_{j=1}^N}\), respectively.
Since \(K_N\)
is skew-symmetric, we have
\[
(K_NZ_l,Z_l)_{\mathbb C^N}
=\sqrt{-1}\,(K_NY_l,X_l)_{\mathbb R^N}.
\]
For \(U:=PZ=(\beta_1+\sqrt{-1}\widetilde\beta_1)/\sqrt2\),
the symmetry of \(E_C\) and the skew-symmetry of \(D_C\) give
\[
(B_CU,U)_{\mathbb C^n}
=\frac12\bigl\{(E_C\beta_1,\beta_1)_{\mathbb R^n}
+(E_C\widetilde\beta_1,\widetilde\beta_1)_{\mathbb R^n}\bigr\}
-\sqrt{-1}\,(D_C\widetilde\beta_1,\beta_1)_{\mathbb R^n}.
\]
By combining these identities, we obtain
\begin{equation}\label{eq:Sm}
\mathcal{S}_N
=\bigl((-2\Lambda\otimes K_N+P^{*}B_CP)Z,Z\bigr)_{\mathbb C^{nN}}.
\end{equation}
Moreover, the Hermitian part of
\(I_{nN}+2\Lambda\otimes K_N-P^{*}B_CP\) is \(I_{nN}-P^{*}E_CP\), and it is
positive definite because \(PP^{*}=vI_n\) and \(I_n-vE_C>0\). Hence
we may apply Lemma~\ref{lem:complex-gauss} with \(q=nN\) to obtain
\begin{equation}\label{eq:expectation-SN}
\mathbb{E}[e^{\mathcal{S}_N}]
=\det(I_{nN}+2\Lambda\otimes K_N-P^{*}B_CP)^{-1}.
\end{equation}

To compute the determinant in \eqref{eq:expectation-SN}, we put
\(
Q_N:=I_{nN}+2\Lambda\otimes K_N
\).
Since \(\Lambda\) is diagonal, \(Q_N\) is the block diagonal matrix
\begin{equation}\label{eq:QN-block}
Q_N=\operatorname{diag}(I_N+2\lambda_lK_N)_{l=1}^{n}.
\end{equation}
For \(a\in\mathbb R\), we define
\begin{align*}
\Delta_N(a)&:=\det(I_N+a^2K_N^*K_N)^{1/2}, \quad
m_N(a):=v\bigl(
(I_N+a^2K_N^*K_N)^{-1}\mathbf{e}_1,\mathbf{e}_1
\bigr)_{\mathbb{R}^N},
\end{align*}
and \(M_N:=\operatorname{diag}(m_N(2\lambda_l))_{l=1}^{n}\).
Since \(e\) is represented by \(\sqrt{v}\,\mathbf{e}_1\) in
\(\mathcal H_N\), these are the counterparts of \(\Delta_R\), \(m_R\)
and \(M\) for \(K_N\). We first show that
\begin{equation}\label{eq:QN-data}
\det Q_N=\prod_{l=1}^n\Delta_N(2\lambda_l),\qquad
PQ_N^{-1}P^{*}=M_N.
\end{equation}
Since \(K_N\) is real and skew-symmetric, we have
\(K_N^{*}K_N=-K_N^{2}\), and hence
\begin{equation}\label{eq:KN-factor}
(I_N-aK_N)(I_N+aK_N)=I_N+a^{2}K_N^{*}K_N
\end{equation}
for \(a\in\mathbb R\). The right-hand side is positive definite, so
that \(I_N+aK_N\) is nonsingular. By taking the determinant of
\eqref{eq:KN-factor} and using
\[\det(I_N-aK_N)=\det\tp \,{(I_N+aK_N)}=\det(I_N+aK_N),\]
we have
\[
\det(I_N+aK_N)^{2}=\Delta_N(a)^{2}.
\]
Since \(\det(I_N+aK_N)\) is continuous in \(a\), does not vanish and
is equal to \(1\) at \(a=0\), it is positive. Hence we obtain
\(\det(I_N+aK_N)=\Delta_N(a)\). Combining this with
\eqref{eq:QN-block}, we have the first identity in \eqref{eq:QN-data}.
Next, \eqref{eq:KN-factor} implies
\[
(I_N+aK_N)^{-1}
=(I_N+a^{2}K_N^{*}K_N)^{-1}
-aK_N(I_N+a^{2}K_N^{*}K_N)^{-1}.
\]
Here the matrix \(K_N(I_N+a^{2}K_N^{*}K_N)^{-1}\) is skew-symmetric,
because \(K_N\) is skew-symmetric and commutes with the symmetric
matrix \((I_N+a^{2}K_N^{*}K_N)^{-1}\). Hence its \((1,1)\)-component
vanishes, and we have
\[
v\bigl((I_N+aK_N)^{-1}\mathbf{e}_1,\mathbf{e}_1\bigr)_{\mathbb R^N}
=m_N(a).
\]
Since \(P=\sqrt{v}\,I_n\otimes\tp \,{\mathbf{e}_1}\), it follows from
\eqref{eq:QN-block} that
\[
PQ_N^{-1}P^{*}
=\operatorname{diag}\bigl(
v\bigl((I_N+2\lambda_lK_N)^{-1}\mathbf{e}_1,\mathbf{e}_1
\bigr)_{\mathbb R^N}\bigr)_{l=1}^{n},
\]
which is the second identity in \eqref{eq:QN-data}.

We now compute the determinant in \eqref{eq:expectation-SN}. We recall
Sylvester's determinant identity
\[
\det(I_{nN}-FG)=\det(I_n-GF)
\]
for \(F\in\mathbb C^{nN\times n}\) and \(G\in\mathbb C^{n\times nN}\),
which follows from \cite[Theorem 1.3.22]{HJ} by evaluating the characteristic polynomials at $1$
(see also \cite[Problem 1.3.P28]{HJ}). Since \(Q_N\) is nonsingular, by
applying this identity with \(F=Q_N^{-1}P^{*}\) and \(G=B_CP\), and
then using \eqref{eq:QN-data}, we have
\begin{align*}
\det(Q_N-P^{*}B_CP)
&=\det Q_N\cdot\det(I_{nN}-Q_N^{-1}P^{*}B_CP)\\
&=\det Q_N\cdot\det(I_n-B_CPQ_N^{-1}P^{*})\\
&=\prod_{l=1}^n\Delta_N(2\lambda_l)\cdot\det(I_n-B_CM_N).
\end{align*}
Moreover, since
\(B_C=\Lambda^{1/2}\,\tp \,{C}\Lambda^{1/2}\) and the matrices \(\Lambda\)
and \(M_N\) are diagonal, Sylvester's identity and the invariance of
the determinant under transposition give
\begin{align*}
\det(I_n-B_CM_N)
&=\det(I_n-\tp \,{C}\Lambda M_N)\\
&=\det(I_n-\Lambda M_NC)=\det(I_n-C\Lambda M_N).
\end{align*}
Therefore \eqref{eq:expectation-SN} implies
\begin{equation}\label{eq:finite-N}
\mathbb{E}[e^{\mathcal{S}_N}]
=\frac{\prod_{l=1}^n\Delta_N(2\lambda_l)^{-1}}{\det(I_n-C\Lambda M_N)}.
\end{equation}

Next, we show that the denominator in \eqref{eq:main} is strictly
positive. By Sylvester's identity, it is equal to
\(\det\{I_n-M^{1/2}B_CM^{1/2}\}\).
Let \(\kappa_C:=\max\{0,\lambda_{\max}(E_C)\}\), where
\(\lambda_{\max}(E_C)\) denotes the largest eigenvalue of \(E_C\).
Since \(0\leq M\leq vI_n\) and \(I_n-vE_C>0\), we have
\[
I_n-M^{1/2}E_CM^{1/2}\geq(1-v\kappa_C)I_n>0.
\]
This is the symmetric part of the real matrix
\(I_n-M^{1/2}B_CM^{1/2}\).
Hence every real eigenvalue of this matrix
is positive, and its nonreal eigenvalues appear in complex conjugate
pairs. Therefore its determinant is strictly positive.

Finally, we let \(N\to\infty\). Since \(K_N\to K\) in the
Hilbert--Schmidt norm, we have \(K_N^{*}K_N\to K^{*}K\) in the trace
norm. 
Hence, by the continuity of the Fredholm determinant and the
resolvent, we have \(\Delta_N(a)\to\Delta_R(a)\) and
\(m_N(a)\to m_R(a)\). 
On the other hand, by the definition of \(\mathcal{S}_N\) and
\(\mathcal{A}^{l}(1)=-2S^{l}\), we have
\[
\mathcal{S}_N-\mathcal{S}^{R}_{C,\Lambda}
=-2\sqrt{-1}\sum_{l=1}^{n}\lambda_l(S^{l}_{N}-S^{l}).
\]
Since \(S^{l}_{N}\) is the partial sum of the series in
\eqref{eq:area-series}, Lemma~\ref{lem:area-operator} implies that
\(S^{l}_{N}\to S^{l}\) in \(L^{2}(\mathbb P)\) as \(N\to\infty\). Hence
\(\mathcal{S}_N\) converges to \(\mathcal{S}^{R}_{C,\Lambda}\) in
\(L^{2}(\mathbb P)\), and in particular in probability. By the
continuity of the exponential function, \(e^{\mathcal{S}_N}\) also
converges to \(e^{\mathcal{S}^{R}_{C,\Lambda}}\) in probability.
Moreover, since the terms containing \(\beta_1\) and
\(\widetilde{\beta}_1\) do not depend on \(N\), we have
\[
|e^{\mathcal{S}_N}|
=\exp \Big [ \frac12\big \{ (E_C\beta_1,\beta_1)_{\mathbb R^n}
+(E_C\widetilde\beta_1,\widetilde\beta_1)_{\mathbb R^n}\big \} \Big ]
\]
for every \(N\in\mathbb N\), and the right-hand side is integrable.
By the dominated convergence theorem, we see that
 \(e^{\mathcal{S}_N}\) converges to \(e^{\mathcal{S}^{R}_{C,\Lambda}}\)
in \(L^{1}(\mathbb P)\). Then by letting
\(N\to\infty\) in \eqref{eq:finite-N}, we obtain \eqref{eq:main}.
This completes the proof.
\end{proof}
%%%%%%%%%%%%%%%%%%%%%%%%%%%%%%%%%%%%%%%%%%%%%%%%%%%%%
\section{Examples}\label{sec:examples}
%%%%%%%%%%%%%%%%%%%%%%%%%%%%%%%%%%%%%%%%%%%%%%%%%%%%%
In this section, we discuss several examples.
\begin{exm}[Fractional Brownian rough path]\label{exam:fbm}
Let \(H\in(1/4,1)\), and let
\[
R_H(s,t)=\tfrac{1}{2}\{s^{2H}+t^{2H}-|t-s|^{2H}\}
\]
be the covariance of fractional Brownian motion with Hurst parameter
\(H\). By \cite{FV, FV11}, \(R_H\) has finite
2D \(\rho\)-variation with
\(\rho=\max\{1,(2H)^{-1}\}<2\), and we have \(R_H(1,1)=1\). Hence our
framework covers the case where the components are independent
fractional Brownian motions. 
We note that L\'{e}vy's 
stochastic areas coincide with those defined by Young integrals when
\(H>1/2\). We write \(\Delta_H=\Delta_{R_H}\), \(m_H=m_{R_H}\) and
\(\mathcal{S}^{H}_{C,\Lambda}=\mathcal{S}^{R_H}_{C,\Lambda}\). Then
Theorem~\ref{thm:main} implies that, under \(I_n-E_C>0\),
\begin{equation}\label{eq:fbm}
\mathbb{E}[e^{\mathcal{S}^{H}_{C,\Lambda}}]
=\frac{\prod_{l=1}^{n}\Delta_H(2\lambda_l)^{-1}}
{\det(I_n-C\Lambda M_H(\Lambda))},\qquad
M_H(\Lambda)=\operatorname{diag}(m_H(2\lambda_l))_{l=1}^{n}.
\end{equation}
\end{exm}
%%%%%%%%%%%%%%%%%%%%%%%%%%%%%%%%%%%%%%%
\begin{exm}[Brownian motion and the AAFN formula]\label{exam:bm}
Let \(R_{\mathrm{BM}}(s,t)=s\wedge t\). Then the components are
independent standard Brownian motions. This is the case \(H=1/2\) in
Example~\ref{exam:fbm}, and we identify \(\mathcal{H}^{R_{\mathrm{BM}}}\) 
with \(L^{2}([0,1])\).
Then we have \(e=1\). As in the proof of
Lemma~\ref{lem:area-operator}, we obtain
$K_{\mathrm{BM}}:=K=\tfrac{1}{2}(V^{*}-V)$,
where
\begin{equation}\label{eq:op-V}
(Vf)(t)=\int_{0}^{t}f(s)\,ds, \quad 
(V^{*}f)(t)=\int_{t}^{1} f(s)\, ds.
\end{equation}
Hence we have 
$$(Kf)'(t)=-f(t), \quad (Kf)(1)=-(Kf)(0).$$
This implies that the nonzero eigenvalues of \(K\)
are \(\pm\sqrt{-1}/((2k-1)\pi)\), \(k\geq 1\), and hence we obtain
\[
\Delta_{1/2}(a)=\prod_{k=1}^{\infty}\Big[1+a^{2}\Big \{\frac{1}{(2k-1)\pi} \Big \}^{2} \Big]
=\cosh(a/2).
\]

Next, for \(a\neq 0\), we put \(u=(I_{\mathcal H}+a^{2}K^{*}K)^{-1}e\) and \(z=K^{2}u\).
Since \(K^{*}K=-K^{2}\), we have \(u=1+a^{2}z\), and \(z\) satisfies the differential equation
$$z''(t)-a^{2}z(t)=1$$ with
the boundary conditions \(z(1)=-z(0)\) and \(z'(1)=-z'(0)\). Solving this equation, 
we have $$u(t)=\frac{\cosh(a(t-1/2))}{\cosh(a/2)}.$$ By
integrating \(u\) over \([0,1]\), we obtain 
$$m_{1/2}(a)=\frac{\tanh(a/2)}{a/2}, \quad 
m_{1/2}(0)=1.
$$
The value at \(a=0\) follows by continuity. Since
\(M_{1/2}(\Lambda)=\Lambda^{-1}\tanh\Lambda\), \eqref{eq:fbm} reduces to
\[
\frac{
\det(\cosh\Lambda)^{-1}
}{\det(I_n-C\tanh\Lambda)}
=\frac{1}{\det(\cosh\Lambda-C\sinh\Lambda)}.
\]

Hence, under \(I_n-E_C>0\), we obtain
\begin{equation}\label{eq:aafn}
\mathbb{E}[e^{\mathcal{S}^{1/2}_{C,\Lambda}}]
=\det(\cosh\Lambda-C\sinh\Lambda)^{-1}.
\end{equation}
This is the AAFN formula.
More precisely, if we identify \(\beta=W^1\),
\(\widetilde{\beta}=W^2\) and \(C=A\) in the notation of
\cite[Theorem~2]{AAFN}, then their stochastic area is \(\mathcal A^l(1)=-2S^l\) in our notation.
Thus the functional \eqref{eq:functional} coincides with
\(\hat{S}_{C,\Lambda}(\sqrt{-1})\) in \cite[(2.3)]{AAFN}, and
substituting \(\sigma=\sqrt{-1}\) into \cite[(2.6)]{AAFN} gives
\eqref{eq:aafn}. The expression in the statement of
\cite[Theorem~2]{AAFN} is
\(\det(\cosh\Lambda+A\sinh\Lambda)^{-1}\), which instead corresponds
to replacing \(A\) by \(-A\) in the functional. The minus sign in
\eqref{eq:aafn} follows from their definition \cite[(2.3)]{AAFN} and
their calculation \cite[(2.6)]{AAFN}. 

The formula in \cite[Theorem~2]{AAFN} is stated under
smallness assumptions on the parameters. In contrast, formula~\eqref{eq:aafn}
holds on the full domain \(I_n-E_C>0\), which is also necessary
for absolute integrability.
\end{exm}
%%%%%%%%%%%%%%%%%%%%%%%%%%%%%%%%%%%%%%%%
%%%%%%%%%%%%%%%%%%%%%%%%%%%%%%%%%%%%%%%%
%%%%%%%%%%%%%%%%%%%%%%%%%%%%%%%%%%%%%%%%%
\begin{exm}[Brownian bridge]\label{exam:bb}
Let \(w=\tp \,(w^{1}_t,\ldots, w^{2n}_{t})_{0\leq t\leq1}\)
be a \(2n\)-dimensional standard Brownian bridge given by
\(w_t^j=B_t^j-tB_1^j, \, 1\leq j\leq2n\), where \(B^1,\ldots,B^{2n}\) are independent
standard Brownian motions. The common covariance of the
components is
\[
R_{\mathrm{BB}}(s,t)=s\wedge t-st,
\]
which has finite 2D \(1\)-variation, and we have \(v=0\).
For a step function \(f\) on \([0,1]\), we put
\[
\bar f:=\int_0^1f(t)\,dt,\qquad \Pi f:=f-\bar f.
\]
Then, for step functions \(f\) and \(g\), we have
\[
(f,g)_{\mathcal H^{R_{\mathrm{BB}}}}
=\int_0^1(\Pi f)(t)(\Pi g)(t)\,dt.
\]
Hence the map \(f\mapsto\Pi f\) extends to an isometric isomorphism
from \(\mathcal H^{R_{\mathrm{BB}}}\) onto
\[
L^2_0:=\left\{f\in L^2([0,1])\, ; \int_0^1f(t)\,dt=0\right\},
\]
and we identify \(\mathcal H^{R_{\mathrm{BB}}}\) with \(L^2_0\)
through this map. Note that \(\Pi\) is the orthogonal projection of
\(L^{2}([0,1])\) onto \(L^2_0\), and that
\(e=\Pi\mathbf1_{[0,1]}=0\) under this identification.

Let \(L_r^{\mathrm{BB}}\) and \(L_r^{\mathrm{BM}}\) be the operators
\(L_r\) in the proof of Lemma~\ref{lem:area-operator} for
\(R_{\mathrm{BB}}\) and \(R_{\mathrm{BM}}\), respectively. Since each
step function \(f\) is represented by \(\Pi f\), we have
\(L_r^{\mathrm{BB}}=\Pi L_r^{\mathrm{BM}}\big|_{L^2_0}\). By letting
\(r\to\infty\), we obtain
\[
K_{\mathrm{BB}}:=K
=\Pi K_{\mathrm{BM}}\big|_{L^2_0}
=-\Pi V\big|_{L^2_0},
\]
where
\(V\) is defined by \eqref{eq:op-V}
and
the last equality follows from \(V^*f=-Vf\) for \(f\in L^2_0\).
Hence we have
\[
(K_{\mathrm{BB}}f)'(t)=-f(t),\qquad
(K_{\mathrm{BB}}f)(1)=(K_{\mathrm{BB}}f)(0).
\]
This implies that, on the complexification of \(L^2_0\), the
functions \(e^{2\pi\sqrt{-1}kt}\), \(k\in\mathbb Z\setminus\{0\}\), form
a complete orthonormal basis consisting of eigenvectors of
\(K_{\mathrm{BB}}\) with eigenvalues \(\sqrt{-1}/(2\pi k)\). Hence the
infinite product formula for the hyperbolic sine gives
\[
\Delta_{R_{\mathrm{BB}}}(a)
=\prod_{k=1}^\infty\left(1+\frac{a^2}{4\pi^2k^2}\right)
=\frac{\sinh(a/2)}{a/2},
\]
where
the quotient is defined to be \(1\) at \(a=0\).
Since $e=0$, we also have \(m_{R_{\mathrm{BB}}}=0\).

Consequently, for every real \(n\times n\) matrix \(C\),
Theorem~\ref{thm:main} yields
\begin{equation}\label{eq:bb-formula}
\mathbb E[e^{\mathcal S^{R_\mathrm{BB}}_{C,\Lambda}}]
=\prod_{l=1}^n\frac{\lambda_l}{\sinh\lambda_l}
=\det(\Lambda^{-1}\sinh\Lambda)^{-1}.
\end{equation}
Here the condition \(I_n-vE_C>0\) is always satisfied because
\(v=0\). Moreover, since \(\beta_1=\widetilde{\beta}_1=0\)
\(\mathbb P\)-almost surely, the terms containing \(C\) in
\eqref{eq:functional} vanish. Hence the right-hand side of
\eqref{eq:bb-formula} does not depend on \(C\), and
\(\lvert e^{\mathcal S^{R_\mathrm{BB}}_{C,\Lambda}}\rvert=1\)
\(\mathbb P\)-almost surely. Formula \eqref{eq:bb-formula} is the
Brownian bridge version
of L\'evy's stochastic area formula \cite{Le51} in our setting.
We also note that the stochastic area \(S^{l}\) is nontrivial
although \(v=0\). Indeed, \eqref{eq:HS} implies
\[
\mathbb E[(S^l)^2]
=\|K_{\mathrm{BB}}\|_{\mathrm{HS}(L^{2}_{0})}^2
=2\sum_{k=1}^\infty\frac1{4\pi^2k^2}
=\frac1{12}.
\]
\end{exm}
%%%%%%%%%%%%%%%%%%%%%%%%%%%%%%%%%%%%%%%%%%%%%%%%%%%%
\noindent
{\bf Declaration on the use of AI.} \
AI tools were used to assist in refining the English writing and in improving the clarity and readability of 
the mathematical exposition in the manuscript.
All mathematical arguments and proofs were checked and verified by the author, who takes full responsibility 
for the content of the manuscript.
\vspace{2mm} \\
\noindent
{\bf Acknowledgements.} \ 
The author was partially supported by 
JSPS Grant-in-Aid for Scientific Research (C) No. 26K06833.
%%%%%%%%%%%%%%%%%%%%%%%%%%%%%%%%%%%%%%%%%%%%%%%%%%%%%


\begin{thebibliography}{99}
%
\bibitem[AAFN13]{AAFN}
H.~Aihara, J.~Akahori, H.~Fujii and Y.~Nitta:
\emph{Tau functions of KP solitons realized in Wiener space},
Bull. London Math. Soc. \textbf{45} (2013), pp. 1301--1309.

\bibitem[CQ02]{CQ}
L.~Coutin and Z.~Qian:
\emph{Stochastic analysis, rough path analysis and fractional Brownian motions},
Probab. Theory Related Fields \textbf{122} (2002), pp. 108--140.

\bibitem[FCU11]{FCU}
A.~Ferreiro-Castilla and F.~Utzet:
\emph{L\'evy area for Gaussian processes: A double Wiener--It\^o
integral approach},
Statist. Probab. Lett. \textbf{81} (2011), pp. 1380--1391.

\bibitem[FR14]{FR}
P.~Friz and S.~Riedel:
\emph{Convergence rates for the full Gaussian rough paths},
Ann. Inst. Henri Poincar\'e Probab. Statist.
\textbf{50} (2014), pp. 154--194.

\bibitem[FV10a]{FV10}
P. Friz and N. Victoir:
Multidimensional Stochastic Processes as Rough Paths: Theory and Applications,
Cambridge Studies in Advanced Mathematics {\bf{120}}, Cambridge University Press, 2010.

\bibitem[FV10b]{FV}
P. Friz and N. Victoir:
\emph{Differential equations driven by Gaussian signals},
Ann. Inst. Henri Poincar\'e Probab. Statist.
\textbf{46} (2010), pp. 369--413.

\bibitem[FV11]{FV11}
P. Friz and N. Victoir:
\emph{A note on higher dimensional \(p\)-variation},
Electron. J. Probab. \textbf{16} (2011), pp. 1880--1899.

\bibitem[Gav77]{Ga}
B. Gaveau:
\emph{Principe de moindre action, propagation de la chaleur et
estim\'ees sous elliptiques sur certains groupes nilpotents},
Acta Math. \textbf{139} (1977), pp. 95--153.

\bibitem[GQ17]{GQ}
X. Geng and Z. Qian:
\emph{Finite dimensional characteristic functions of Brownian rough path},
Front. Math. China \textbf{12} (2017), pp. 859--877.

\bibitem[HS83]{HS}
K. Helmes and A. Schwane:
\emph{L\'evy's stochastic area formula in higher dimensions},
J. Funct. Anal. \textbf{54} (1983), pp. 177--192.

\bibitem[HJ13]{HJ}
R.A. Horn and C.R. Johnson:
Matrix Analysis, 2nd ed.,
Cambridge University Press, 2013.

\bibitem[IKM94]{IKM}
N. Ikeda, S. Kusuoka and S. Manabe:
\emph{L\'evy's stochastic area formula for Gaussian processes},
Comm. Pure Appl. Math. \textbf{47} (1994), pp. 329--360.

\bibitem[IT04]{IT}
N. Ikeda and S. Taniguchi:
\emph{Quadratic Wiener functionals, Kalman--Bucy filters, and the KdV
equation},
in: Stochastic analysis and related topics in Kyoto,
Adv. Stud. Pure Math. \textbf{41}, Math. Soc. Japan, Tokyo, 2004,
pp. 167--187.

\bibitem[L\'ev40]{Le40}
P. L\'evy:
\emph{Le mouvement brownien plan},
Amer. J. Math. \textbf{62} (1940), pp. 487--550.

\bibitem[L\'ev51]{Le51}
P. L\'evy:
\emph{Wiener's random function, and other Laplacian random functions},
in: Proceedings of the Second Berkeley Symposium on Mathematical
Statistics and Probability, 1950, Univ. California Press, Berkeley,
1951, pp. 171--187.

\bibitem[Yor80]{Yo}
M. Yor:
\emph{Remarques sur une formule de Paul L\'evy},
in: S\'eminaire de Probabilit\'es XIV, Lecture Notes in Math.
\textbf{784}, Springer, Berlin, 1980, pp. 343--346.
%
\end{thebibliography}
\end{document}